\documentclass[12pt]{amsart}
\usepackage{lineno,hyperref}

\usepackage{pgf,tikz}
\usetikzlibrary{arrows}
\usepackage{cite}
\usepackage{blindtext}
\usepackage{amsmath, amssymb,graphicx,enumerate,color}
\usepackage{multirow}
\usepackage{multicol}
\usepackage{xcolor}

\def\pmatrix{\left(\begin{matrix}}
\def\endpmatrix{\end{matrix}\right)}

\def\Z{{\mathbb Z}}
\def\C{{\mathbb C}}

\def\det{\operatorname{det}}

\def\H{{\mathcal H}}

\theoremstyle{plain}
\newtheorem{theorem}{Theorem}
\newtheorem{lemma}[theorem]{Lemma}

\theoremstyle{definition}

\allowdisplaybreaks

\begin{document}


\title{ A character of the Siegel modular group of level 2 from theta constants }

\author{Xinhua Xiong}
\address{Department of Mathematics, China Three Gorges University\\
 Yichang, Hubei Province, 443002, P.R. China\\}
 
\email{xinhuaxiong@yahoo.com}

\begin{abstract}  
Given a characteristic, we define a character of the Siegel modular group of level 2, the computations of their values are also obtained. By using our theorems, some key theorems of Igusa [1] can be recovered.
 \end{abstract}

 
 \maketitle

\section{Introduction.}
The theta function of characteristic $m$ of degree $g$ is the series
\begin{eqnarray*}
\theta_{m} (\tau,z):=\sum\limits_{p\in\Z^g} \exp \pi i \left[\left(
p+\frac{m^{\prime}}{2},\tau(p+\frac{m^{\prime}}{2})\right)\right.\\
+2\left.\left(p+\frac{m^{\prime}}{2},z+
\frac{m^{\prime\prime}}{2}\right)\right],
\end{eqnarray*}
where $\tau \in \H_{g}$, $z\in \C^g$, $m=\pmatrix m^{\prime}\\ m^{\prime\prime}\endpmatrix \in \Z^{2g}$, $\H_g$ is the  Siegel upper half-plane, $m^{\prime}$ and $m^{\prime\prime}$ denote vectors in $\Z^{g}$ determined by the first and last $g$ coefficients of $m$. If we put $z=0$, we get the theta constant $\theta_{m}(\tau)=\theta_{m}(\tau,0)$. The study of theta functions and theta constants has a long history, and they are very important objects in arithmetic and geometry. They can be used to construct modular forms and to study geometric properties of abelian varieties. Farkas and Kra's book [1] contains very detailed descriptions for the case of degree one. In [2], [3], and [4], Matsuda gives new formulas and applications. It is Igusa in [5] who began to study the cases of higher degrees. He used $\theta_{m}(\tau)\theta_{n}(\tau)$ to determine the structure of the graded rings of modular forms belonging to the group $\Gamma_{g}(4,8)$. 

In this note, we will  define a character of the group $\Gamma_{g}(2)$, the principal congruence group of degree $g$ and of level $2$. We obtained its computation formula.   Using our results, Igusa's key Theorem 3 in [5] can be recovered.

 \section{The Siegel modular group of level $2$.}
The Siegel modular group ${\rm Sp}(g,\Z)$ of degree $g$ is the group of $2g \times 2g$ integral matrices $M$ satisfying
$$
M\pmatrix 0& I_{g}\\ -I_{g}&0 \endpmatrix   {^{t}M}=\pmatrix 0& I_{g}\\ -I_{g}&0 \endpmatrix,
$$
in which ${^{t}M}$ is the transposition of $M$, $I_{g}$ is the identity of degree $g$. If we put $M=\pmatrix a& b\\ c&d \endpmatrix$, the condition for $M$ in ${\rm Sp}(g, \Z)$ is $a^{t}d-b^{t}c=I_{g}$, $a^{t}b$ and ${c^{t}d}$ are symmetric matrices. In fact, if $M$ is in ${\rm Sp}(g, \Z)$, then $^{t}bd$ and ${^{t}ac}$ are also symmetric, see [6, p. 437]. In this paper, we discuss two special subgroups of the Siegel modular group. The first is the principal congruence subgroup $\Gamma_{g}(2)$ of degree $g$ and of level $2$ which is defined by $M\equiv I_{2g} \pmod{2}$. The second is the Igusa modular group $\Gamma_{g}(4,8)$, which is defined by $M\equiv I_{2g}\pmod{4}$ and $(a^{t}b)_{0}\equiv (c^{t}d)_{0}\equiv 0 \pmod{8}$. Where if $s$ is a square matrix, we arrange its diagonal coefficients in a natural order to form a vector $(s)_{0}$.  The Siegel modular group ${\rm Sp}(g,\Z)$ acts on $\H_{g}$ by the formula
$$
M\tau=(a\tau+b)(c\tau +d)^{-1},
$$
where $ \tau \in \H_{g}, M=\pmatrix a& b\\ c&d \endpmatrix \in {\rm Sp}(g, \Z)$.
An element $m\in \Z^{2g}$ is called a theta characteristic of degree $g$. If $n$ is another characteristic, then we have 
\begin{equation*}
\theta_{m+2n}(\tau, z)=(-1)^{^t{m^{\prime}n^{\prime\prime}}}\theta_{m}(\tau, z).
\end{equation*}
Since $$\theta_{m}(\tau, -z)=(-1)^{^t{m^{\prime}m^{\prime\prime}}}\theta_{m}(\tau, z),$$
$m$ is called even or odd according as ${^{t}m^{\prime}m^{\prime\prime}}$ is even or odd. Only for even $m$, theta constants are none zero. Given a characteristic $m$ and an element $M$ in the Siegel modular group, we define 
$$
M\circ m=
\pmatrix d &  -c \\
-b  & a \endpmatrix \pmatrix m^{\prime}\\ m^{\prime\prime} \endpmatrix +
\pmatrix (c^{t}d)_{0}\cr (a^{t}b)_{0} \endpmatrix.
$$
This operation modulo $2$ is a group action, i.e. $M_{1}\circ(M_{2}\circ m) \equiv (M_{1}M_{2})\circ m \pmod{2}$. Next we explain the transformation formulas for theta constants: for any $m \in \Z^{2g}$ and $M \in {\rm Sp}(g, \Z)$, we put
\begin{eqnarray*}
\Phi_{m}(M):=-\frac{1}{8}\left( ^{t}m^{\prime} {^{t}bd}m^{\prime} + {^{t}m^{\prime\prime}} {^{t}ac} m^{\prime\prime}\right.\\
\left.-2 {^{t}m^{\prime}}{^{t}bc}m^{\prime\prime}-2{^{t}(a^{t}b)_{0}}(dm^{\prime}-cm^{\prime\prime})
\right),
\end{eqnarray*}
then we have 
\begin{eqnarray*}
\theta_{M\circ m}(M\tau)=a(M)e(\Phi_{m}(M))\det(c\tau+d)^{\frac{1}{2}}\theta_{m}(\tau),
\end{eqnarray*}
in which $a(M)$ is an eighth root of unity depending only on $M$ and the choice of square root sign for $\det(c\tau+d)^{\frac{1}{2}}$, and
$e(\Phi_{m}(M))=e^{2\pi i \Phi_{m}(M)}$. From now on, we always discuss the group $\Gamma_{g}(2)$, unless specified. Hence, we can write
$$
\Phi_{m}(M)=-\frac{1}{8}\left( ^{t}m^{\prime} {^{t}bd}m^{\prime} + {^{t}m^{\prime\prime}} {^{t}ac} m^{\prime\prime}-2{^{t}(a^{t}b)_{0}}dm^{\prime}\right).
$$
\section{Main theorems and proofs.}
 The character mentioned in the abstract is as follows:

\noindent{\bf Definition. }
 Let $m\in \Z^{2g}, M \in \Gamma_{g}(2)$, we define $\chi_{m}(M)$ by 
\begin{equation*}
\theta_{m}(M\tau)=a(M)\chi_{m}(M)\det(c\tau +d)^{\frac{1}{2}}\theta_{m}(\tau),
\end{equation*}
where $a(M)$ comes from the transformation formulas of theta constants.

\begin{theorem}
For a fixed $m$, $\chi_{m}(M)$ is a character of $\Gamma_{g}(2)$.
\end{theorem}
The proof of Theorem 3.1 needs two lemmas, in which $\Gamma_{g}(2)$ is essential.
\begin{lemma}
If $m, n \in \Z^{2g}$ and $m\equiv n \pmod{2}$, then $\Phi_{m}(M)\equiv \Phi_{n}(M) \pmod{1}$.
\end{lemma}
\noindent{\it Proof.} Let $m=n+2\Delta$, then 
\begin{eqnarray*}
&&^{t}m^{\prime} {^{t}bd}m^{\prime}\\
&=&{^{t}(n^{\prime}+2\Delta^{\prime})} {^{t}bd}(n^{\prime}+2\Delta^{\prime})\\
&=&^{t}n^{\prime} {^{t}bd}n^{\prime}+4^{t}\Delta^{\prime} {^{t}bd}\Delta^{\prime}+2^{t}n^{\prime} {^{t}bd}\Delta^{\prime}+2^{t}\Delta^{\prime} {^{t}bd}n^{\prime}\\
&\equiv& ^{t}n^{\prime} {^{t}bd}n^{\prime} \pmod{8}, 
\end{eqnarray*}
since $b\equiv 0 \pmod{2}$ and $^{t}bd$ is symmetric, the last two terms are equal. Similarly, $^{t}ac$ is symmetric which implies that 
$^{t}m^{\prime\prime} {^{t}ac}m^{\prime\prime}
\equiv {^{t}n^{\prime\prime}} {^{t}ac}n^{\prime\prime} \pmod{8} $.  Moreover, $2{^{t}(a^{t}b)_{0}}dm^{\prime}\equiv 2{^{t}(a^{t}b)_{0}}dn^{\prime}\pmod{8}$ is trivial. By the definition of $\Phi_{m}(M)$, Lemma 3.2 is true.
\begin{lemma}
For $M, M^{\prime} $ $\in $ $\Gamma_{g}(2)$, we have $\Phi_{M^{\prime}\circ m}(M)\equiv \Phi_{m}(M) \pmod{1}$.
\end{lemma}
\noindent {\it Proof.}  This lemma can be proved from the definition of $M^{\prime}\circ m$, $M^{\prime}$ is in $\Gamma_{g}(2)$ and Lemma 3.2. 

\noindent{\it Proof of Theorem 3.1.} 

We firstly give a formula for $\chi_{m}(M)$. By the definition of the operation $\circ $, we can find a unique
$n$ in $\Z^{2g}$ with $M\circ n=m$. Define $\Delta\in \Z^{2g}$ by $m+2\Delta=n$, then by Lemma 3.2 and Lemma 3.3, we have
\begin{eqnarray*}
&&\theta_{m}(M\tau)\\
&=& \theta_{M\circ n}(M\tau)\\
                            &=& a(M)e(\Phi_{n}(M))\det(c\tau +d)^{\frac{1}{2}}\theta_{n}(\tau) \\
                           &=& a(M)e(\Phi_{m}(M))\det(c\tau +d)^{\frac{1}{2}}\theta_{n}(\tau) \\
                           &=& a(M)e(\Phi_{m}(M))\det(c\tau +d)^{\frac{1}{2}}\theta_{m+2\Delta}(\tau)   \quad\quad\\
                           &=& a(M)e(\Phi_{m}(M))\det(c\tau +d)^{\frac{1}{2}} (-1)^{^{t}m^{\prime}\Delta^{\prime\prime}}\theta_{m}(\tau). 
\end{eqnarray*}
Hence, 
\begin{equation}
\chi_{m}(M)=e(\Phi_{m}(M))(-1)^{^{t}m^{\prime}\Delta^{\prime\prime}}.
\end{equation}

To prove Theorem 3.1 is equivalent to prove $\chi_{m}(M_{1})\chi_{m}(M_{2})=\chi_{m}(M_{1}M_{2})$ for any $M_{1}, M_{2} \in \Gamma_{g}(2)$. Now fix $M_{1},M_{2}$ and $m$, define $\Delta_{1}, \Delta_{2}$ by $m+2\Delta_{1}=n_{1}, M_{1}\circ n_{1}=m; \,m+2\Delta_{2}=n_{2}, M_{2}\circ n_{2}=m.$  Write 
$$
M_{1}=  \pmatrix a_{1} &  b_{1} \\ c_{1}  & d_{1} \endpmatrix,   M_{2}=\pmatrix a_{2} &  b_{2} \\ c_{2}  & d_{2} \endpmatrix,  
$$
by (1), we have $$\chi_{m}(M_{1})=e(\Phi_{m}(M_{1}))(-1)^{^{t}m^{\prime}\Delta_{1}^{\prime\prime}}$$ and $$\chi_{m}(M_{2})=e(\Phi_{m}(M_{2}))(-1)^{^{t}m^{\prime}\Delta_{2}^{\prime\prime}}.$$  In order to compute $\chi_{m}(M_{1}M_{2})$,  we write $ \tau^{\prime}=M_{2}\tau,\,M=M_{1}M_{2}=\pmatrix a &  b \\  c & d \endpmatrix,$ and define $\Delta_{3}$ by $m+2\Delta_{3}=\bar{n}$ with $M_{1}\circ n_{1}=m, M_{2}\circ \bar{n}=n_{1}$. Then by Lemma 3.2 and 3.3, we have
\begin{eqnarray*}
&&\theta_{m}(M_{1}M_{2}\tau) \\
&=& \theta_{M_{1}\circ n_{1}}(M_{1}\tau^{\prime})\\
 &=& a(M_{1})e(\Phi_{n_{1}}(M_{1}))\det(c_{1}\tau^{\prime} +d_{1})^{\frac{1}{2}}\theta_{n_{1}}(\tau^{\prime})\\ 
 &=& a(M_{1})e(\Phi_{m}(M_{1}))\det(c_{1}\tau^{\prime} +d_{1})^{\frac{1}{2}}\theta_{M_{2}\circ \bar{n}}(M_{2}\tau)\\ 
 &=& a(M_{1})e(\Phi_{m}(M_{1}))\det(c_{1}\tau^{\prime} +d_{1})^{\frac{1}{2}} a(M_{2})e(\Phi_{\bar{n}}(M_{2}))\\
 &&\times\det(c_{2}\tau +d_{2})^{\frac{1}{2}} \theta_{\bar{n}}(\tau)\\ 
 &=& a(M_{1})a(M_{2})e(\Phi_{m}(M_{1}))e(\Phi_{n_{1}}(M_{2}))\det(c_{1}\tau^{\prime} +d_{1})^{\frac{1}{2}}\\
  &&\times\det(c_{2}\tau +d_{2})^{\frac{1}{2}}\theta_{\bar{n}}(\tau)\\
  &=& a(M_{1})a(M_{2})e(\Phi_{m}(M_{1}))e(\Phi_{m}(M_{2}))\\
 &&\times \det(c\tau +d)^{\frac{1}{2}} \theta_{m+2\Delta_{3}}(\tau)\\
  &=& a(M_{1}M_{2})e(\Phi_{m}(M_{1}))e(\Phi_{m}(M_{2}))\\
  &&\times\det(c\tau +d)^{\frac{1}{2}} (-1)^{^{t}m^{\prime}\Delta_{3}^{\prime\prime}} \theta_{m}(\tau),
 \end{eqnarray*}
in which, we use $c\tau+d=(c_{1}\tau^{\prime}+d_{1})(c_{2}\tau+d_{2})$ and $a(M)=a(M_{1})a(M_{2})$, the later is implied by Igusa's Theorem 2 in [5]. Hence,
$$
\chi_{m}(M_{1}M_{2})=e(\Phi_{m}(M_{1}))e(\Phi_{m}(M_{2}))(-1)^{^{t}m^{\prime}\Delta_{3}^{\prime\prime}}.
$$
Now we compute $$(-1)^{^{t}m^{\prime}\Delta_{1}^{\prime\prime}}, (-1)^{^{t}m^{\prime}\Delta_{2}^{\prime\prime}} \mbox{ and } (-1)^{^{t}m^{\prime}\Delta_{3}^{\prime\prime}}.$$
 From $M_{1}\circ n_{1}=m$, we get
\begin{equation}
n_{1}=\pmatrix  ^{t}a_{1}m^{\prime}+{^{t}c_{1}}m^{\prime\prime}-{^{t}a_{1}(c_{1}}{^{t}d_{1}})_{0}-{^{t}c_{1}}(a_{1}^{t}b_{1})_{0}\\ ^{t}b_{1}m^{\prime}+{^{t}d_{1}}m^{\prime\prime}-{^{t}b_{1}(c_{1}}{^{t}d_{1}})_{0}-{^{t}d_{1}}(a_{1}^{t}b_{1})_{0} \endpmatrix.
\end{equation} Therefore, from $m+2\Delta_{1}=n_{1}$, we have
\begin{eqnarray*}
&&\Delta_{1}^{\prime\prime}
=\frac{n_{1}^{\prime\prime}-m^{\prime\prime}}{2}\\
&=&\frac{^{t}b_{1}m^{\prime}}{2}+\frac{(^{t}d_{1}-I_{g})m^{\prime\prime}}{2}-\frac{{^{t}b_{1}(c_{1}}{^{t}d_{1}})_{0}}{2}-\frac{{^{t}d_{1}}(a_{1}^{t}b_{1})_{0}}{2}\\
&\equiv&  \frac{^{t}b_{1}m^{\prime}}{2}+\frac{(^{t}d_{1}-I_{g})m^{\prime\prime}}{2}-\frac{(^{t}b_{1})_{0}}{2} \pmod{2},
\end{eqnarray*}
by noting that $b_{1}\equiv c_{1}\equiv 0 \pmod{2}$ and $a_{1}\equiv d_{1}\equiv I_{g}\pmod{2}$. So
\begin{eqnarray*}
&&{^{t}m^{\prime}\Delta_{1}^{\prime\prime}}\\
&\equiv&  \frac{^{t}m^{\prime}{^{t}b_{1}}m^{\prime}}{2}+\frac{^{t}m^{\prime}(^{t}d_{1}-I_{g})m^{\prime\prime}}{2}-\frac{^{t}m^{\prime}(^{t}b_{1})_{0}}{2} \pmod{2} \\
&\equiv& \frac{^{t}m^{\prime}(^{t}d_{1}-I_{g})m^{\prime\prime}}{2} \pmod{2},
\end{eqnarray*}
because $\frac{^{t}b_{1}}{2} \pmod{2}$ is symmetric, which follows from the fact $a^{t}b_{1}$ is symmetric. Similarly, $${^{t}m^{\prime}\Delta_{2}^{\prime\prime}}\equiv \frac{^{t}m^{\prime}(^{t}d_{2}-I_{g})m^{\prime\prime}}{2} \pmod{2}.$$
The computation for $\Delta_{3}^{\prime\prime}$ is more complicated. Recall that $M_{1}\circ n_{1}=m,\, M_{2}\circ \bar{n}=n_{1}, m+
2\Delta_{3}=\bar{n}$ and (2), we have
\begin{eqnarray*}
n_{1}&\equiv& \pmatrix  ^{t}a_{1}m^{\prime}+{^{t}c_{1}}m^{\prime\prime}-{^{t}a_{1}(c_{1}}{^{t}d_{1}})_{0}\\ ^{t}b_{1}m^{\prime}+{^{t}d_{1}}m^{\prime\prime}-{^{t}d_{1}}(a_{1}^{t}b_{1})_{0} \endpmatrix\pmod{4}.
\end{eqnarray*}
From $M_{2}\circ \bar{n}=n_{1}$, we get
$$
\bar{n} \equiv \pmatrix  ^{t}a_{2}n_{1}^{\prime}+{^{t}c_{2}}n_{1}^{\prime\prime}-{^{t}a_{2}(c_{2}}{^{t}d_{2}})_{0}\\ ^{t}b_{2}n_{1}^{\prime}+{^{t}d_{2}}n_{1}^{\prime\prime}-{^{t}d_{2}}(a_{2}^{t}b_{2})_{0} \endpmatrix\pmod{4}.
$$
Therefore,
\begin{eqnarray*}
\bar{n}^{\prime\prime}
 &\equiv& ^{t}b_{2}(^{t}a_{1}m^{\prime}+{^{t}c_{1}}m^{\prime\prime}-{^{t}a_{1}(c_{1}}{^{t}d_{1}})_{0}) +{^{t}d_{2}}(^{t}b_{1}m^{\prime}\\
&&+{^{t}d_{1}}m^{\prime\prime}-{^{t}d_{1}}(a_{1}^{t}b_{1})_{0})-{^{t}d_{2}}(a_{2}^{t}b_{2})_{0} \pmod{4}\\
&\equiv& {^{t}b_{2}}{^{t}a_{1}}m^{\prime}+{^{t}d_{2}}{^{t}b_{1}}m^{\prime}+{^{t}d_{2}}{^{t}d_{1}}m^{\prime\prime}\\
&&-{^{t}d_{2}}{^{t}d_{1}}(a_{1}^{t}b_{1})_{0}-{^{t}d_{2}}(a_{2}{^{t}b_{2}})_{0}\pmod{4}.
\end{eqnarray*}
By the definition of $\Delta_{3}$, we have
\begin{eqnarray*}
&&\Delta_{3}^{\prime\prime}
=\frac{\bar{n}^{\prime\prime}-m^{\prime\prime}}{2}\\
&\equiv& \frac{{^{t}b_{2}}{^{t}a_{1}}m^{\prime}}{2}+\frac{{^{t}d_{2}}{^{t}b_{1}}m^{\prime}}{2}+\frac{({^{t}d_{2}}{^{t}d_{1}}-I_{g})m^{\prime\prime}}{2}\\
&&-\frac{{^{t}d_{2}}{^{t}d_{1}}(a_{1}^{t}b_{1})_{0}}{2}-\frac{{^{t}d_{2}}(a_{2}{^{t}b_{2}})_{0}}{2}\pmod{2}\\
&\equiv& \frac{{^{t}b_{2}}m^{\prime}}{2}+\frac{{^{t}b_{1}}m^{\prime}}{2}+\frac{({^{t}d_{2}}{^{t}d_{1}}-I_{g})m^{\prime\prime}}{2}\\
&&-\frac{(^{t}b_{1})_{0}}{2}-\frac{({^{t}b_{2}})_{0}}{2}\pmod{2},\\
and\\
&&^{t}m^{\prime}\Delta_{3}^{\prime\prime}\\
& \equiv & \frac{^{t}m^{\prime}{^{t}b_{2}}m^{\prime}}{2}+\frac{^{t}m^{\prime}{^{t}b_{1}}m^{\prime}}{2}+\frac{^{t}m^{\prime}({^{t}d_{2}}{^{t}d_{1}}-I_{g})m^{\prime\prime}}{2}\\
 &&-\frac{^{t}m^{\prime}(^{t}b_{1})_{0}}{2}-\frac{^{t}m^{\prime}({^{t}b_{2}})_{0}}{2}\pmod{2}\\
&\equiv& \frac{^{t}m^{\prime}({^{t}d_{2}}{^{t}d_{1}}-I_{g})m^{\prime\prime}}{2} \pmod{2},
\end{eqnarray*}
by using the expansions of quadratic forms and the fact that $\frac{^{t}b_{2}}{2}, \frac{^{t}b_{1}}{2}$ are symmetric modulo $2$.
Finally, the verification of $$^{t}m^{\prime}\Delta^{\prime\prime}\equiv{^{t}m^{\prime}}\Delta_{1}^{\prime\prime}+{^{t}m^{\prime}}\Delta_{2}^{\prime\prime} \pmod{2}$$ is easy, which comes from the simple fact
$$
(^{t}d_{2}-I_{g})(^{t}d_{1}-I_{g}) \equiv 0 \pmod{4}.
$$
This completes the proof of Theorem 3.1.

In [5], Igusa gave the generators  $A_{ij}, B_{ij}, C_{ij}$ of $\Gamma_{g}(2)$, where
\begin{enumerate}
\item $1\leq i \neq j\leq g$, $ A_{ij} =\pmatrix a &  0 \\ 0  & d \endpmatrix, d={^{t}a^{-1}}, a$ is obtained by replacing $(i,j)$-coefficient in $I_{g}$ by $2$;
\item $1\leq i \leq g$, $ A_{ii} =\pmatrix a &  0 \\ 0  & d \endpmatrix, d={^{t}a^{-1}}, a$ is obtained by replacing $(i,i)$-coefficient in $I_{g}$ by $-1$;
\item $1\leq i < j\leq g$, $ B_{ij} =\pmatrix I_{g} &  b \\ 0  & I_{g} \endpmatrix, b$ is obtained by replacing $(i,j)$- and $(j,i)$-coefficients in $0$ by $2$;
\item $1\leq i \leq g$, $ B_{ii} =\pmatrix I_{g} &  b \\ 0  & I_{g} \endpmatrix, b$ is obtained by replacing $(i,i)$-coefficient in $0$ by $2$;
\item $1\leq i \leq j\leq g$, $ C_{ij} ={^{t}B_{ij}}$.
\end{enumerate}
By noting the computation of $\Delta$, which depends on $m$ and $M$, we find $(-1)^{^{t}m^{\prime}}\Delta^{\prime\prime}=1$ for $M=B_{ij}$ or $C_{ij}$, because in these cases, $d=I_{g}$, hence ${^{t}m^{\prime}}\Delta^{\prime\prime}=0$. If $M=A_{ij}$, it is easy to find that ${^{t}m^{\prime}}\Delta^{\prime\prime}=-m_{i}^{\prime}m_{j}^{\prime\prime}\equiv m_{i}^{\prime}m_{j}^{\prime\prime} \pmod{2}$. We can easily compute $\Phi_{m}(M)$ for 
$M=A_{ij}, B_{ij}$  and $C_{ij}$. Now  the values of $\chi_{m}(M)$ for the generators are
\begin{eqnarray*}
&\chi_{m}(A_{ij})=(-1)^{m_{i}^{\prime}m_{j}^{\prime\prime}}, \chi_{m}(B_{ij})=(-1)^{m_{i}^{\prime}m_{j}^{\prime}}, \\
&\chi_{m}(B_{ii})=(-1)^{m_{i}^{\prime}}e\left(\frac{-(m_{i}^{\prime})^{2}}{4}\right),\\
& \chi_{m}(C_{ii})=e\left(\frac{-(m_{i}^{\prime\prime})^{2}}{4}\right),
\chi_{m}(C_{ij})=(-1)^{m_{i}^{\prime\prime}m_{j}^{\prime\prime}}.
\end{eqnarray*}
Using the definitions of $\Phi_{m}(M)$ and $^{t}m^{\prime}\Delta^{\prime\prime}$, it is easy to prove $\chi_{m}(M)=1$ for $M$  in the
Igusa modular group $\Gamma_{g}(4,8)$. Hence, by the computations above, we get 
\begin{theorem}
Write $M$ in the form
$$
M=\prod_{1\leq i,j \leq g}A_{ij}^{p_{ij}}\prod_{1\leq i \leq j\leq g}B_{ij}^{q_{ij}}\prod_{1\leq i\leq j \leq g}C_{ij}^{r_{ij}}M^{\prime}
$$
\noindent with $p_{ij}, q_{ij}, r_{ij} \in \Z$ and $M^{\prime}$ is in the commutator subgroup of $\Gamma_{g}(2)$, which is in $\Gamma_{g}(4,8)$, then 
$$
\chi_{m}(M)=(-1)^{A}e(-\frac{B}{4}) 
$$ with
\begin{eqnarray*}
&A&=\sum_{1\leq i,j \leq g}p_{ij}m_{i}^{\prime}m_{j}^{\prime\prime}+\sum_{1\leq i \leq j \leq g}q_{ij}m_{i}^{\prime}m_{j}^{\prime}\\
&&+\sum_{1\leq i < j \leq g}r_{ij}m_{i}^{\prime\prime}m_{j}^{\prime\prime},\\
&B&=\sum_{1\leq i \leq g}q_{ii}(m_{i}^{\prime})^{2}+\sum_{1\leq i \leq g}r_{ii}(m_{i}^{\prime\prime})^{2}.
\end{eqnarray*}
\end{theorem}

\section{Applications.\\}

If we define $\psi(\tau)=\theta_{m}(\tau)\theta_{n}(\tau)$, then for $M \in \Gamma_{g}(2)$,
$$
\psi(M\tau)=a^{2}(M)\chi_{m}(M)\chi_{n}(M)\det(c\tau+d)\theta_{m}(\tau)\theta_{n}(\tau),$$
we find $\chi_{m}(M)\chi_{n}(M)$ is exactly the character defined by Igusa in [5], hence our theorems can recover Igusa's Theorem 3 in [5].  Our character $\chi_{m}(M)$ is more fundamental, moreover we can see the relations between $\chi_{m}(M)$ and $\Phi_{m}(M)$.

We can use our results to give a more transparent proof of the key part of Theorem 5 in Igusa's paper [5]. The key part of Theorem 5 in that paper is  from the invariant condition that $\frac{\theta_{m}(M\tau)}{\theta_{n}(M\tau)}=\frac{\theta_{m}(\tau)}{\theta_{n}(\tau)}$
holds for all even $m, n$ to infer $M$ is in $\Gamma_{g}(4,8)$, here $M$ is in $\Gamma_{g}(2)$.  By the definition of $\chi_{m}(M)$, this is equivalent to the congruence $\chi_{m}(M)=\chi_{n}(M)$ holds for all even $m, n$, i.e.
$$
e(\Phi_{m}(M))(-1)^{\frac{^{t}m^{\prime}(^{t}d_{}-I_{g})m^{\prime\prime}}{2}}$$ is equal to 
$$e(\Phi_{n}(M))(-1)^{\frac{^{t}n^{\prime}(^{t}d_{}-I_{g})n^{\prime\prime}}{2}}.
$$
Let $n^{\prime}=n^{\prime\prime}=0$, we get for any even $m$, 
$$
\frac{^{t}m^{\prime}(^{t}d_{}-I_{g})m^{\prime\prime}}{2}\equiv 0 \pmod{2}
$$
and 
$$
^{t}m^{\prime} {^{t}bd}m^{\prime} + {^{t}m^{\prime\prime}} {^{t}ac} m^{\prime\prime}-2{^{t}(a^{t}b)_{0}}dm^{\prime}\equiv 0\pmod{8}.
$$
The first congruence implies $d\equiv I_{g} \pmod{4}$, from $a^{t}d-b^{t}c=I_{g}$, we get $a\equiv I_{g} \pmod{4}$. In the second congruence, let $m^{\prime}=0$, we get ${^{t}m^{\prime\prime}}{^{t}a}cm^{\prime\prime}\equiv 0 \pmod{8}$, this implies $^{t}ac\equiv 0 \pmod{4}$, hence $c\equiv 0 \pmod{4}$,  and $(^{t}ac)_{0}\equiv 0 \pmod{8}$. If we write ${^{t}a}=I_{g}+4\bar{a}$, then we have 
$$
({^{t}a}c)_{0}=((I_{g}+4\bar{a}) c)_{0}\equiv (c)_{0}\equiv 0\pmod{8}.
$$
Let $m^{\prime\prime}=0$,  we get the congruence $${^{t}m^{\prime}}{^{t}b}dm^{\prime}-2{^{t}(a^{t}b)_{0}}dm^{\prime}\equiv 0 \pmod{8},$$ which is equivalent 
to the congruence
\begin{equation}
{^{t}m^{\prime}}{^{t}b}dm^{\prime}-2{^{t}(b)_{0}}m^{\prime}\equiv 0 \pmod{8}.
\end{equation}
Write ${^{t}b}=(2\overline{b_{ij}})$ and $d=I_{g}+4\bar{d}$, then 
\begin{eqnarray*}
&&{^{t}m^{\prime}}{^{t}b}dm^{\prime}-2{^{t}(b)_{0}}m^{\prime}\\
&=&{^{t}m^{\prime}}{2\overline{b_{ij}}}(I_{g}+4\overline{d})m^{\prime}-2(2\overline{b_{ij}})_{0}m^{\prime}\\
&\equiv& 2{^{t}m^{\prime}}{\overline{b_{ij}}}m^{\prime}-4(\overline{b_{ij}})_{0}m^{\prime}\pmod{8}\\
&\equiv& 2\overline{b_{11}}{m_{1}^{\prime}}^{2}-4\overline{b_{11}}m_{1}^{\prime}\pmod{8}
\end{eqnarray*}
by taking ${^{t}m^{\prime}}=(m_{1}^{\prime},0,0,\,\cdots,0)$. Hence (3) implies $\overline{b_{11}}\equiv 0 \pmod{4}$. Similarly, we can prove
$\overline{b_{ii}}\equiv 0 \pmod{4}$ holds for each $1\leq i \leq g$. Therefore, we have $(b)_{0}\equiv 0\pmod{8}$. Combing it with (3), we find the congruence ${^{t}m^{\prime}}{^{t}b}dm^{\prime}\equiv 0 \pmod{8}$ holds for any even $m^{\prime}$, which implies ${^{t}b}d\equiv 0\pmod{4}$, hence $b\equiv 0\pmod{4}$. The analysis above shows that $M$ is in $\Gamma_{g}(4)$ and $(b)_{0}\equiv (c)_{0}\equiv 0 \pmod{8}$.  The observation of Igusa in [5, p. 222, line 4-line 6] shows $M$ is in $\Gamma_{g}(4,8)$.
 
 \noindent{\bf Acknowledgments.}  The author would like to thank the referee for his/her helpful corrections and suggestions.

 \let\cleardoublepage\clearpage

\end{document}